\documentclass{article} 
\usepackage{amsmath,amssymb,amsfonts}   
\usepackage{graphicx,bm}

\renewcommand{\v}{\widetilde{v}}
\renewcommand{\u}{\widetilde{u}}
\newcommand{\calU}{{\mathcal U}}
\newcommand{\calG}{{\mathcal G}}
\newcommand{\R}{{\mathbb R}}
\renewcommand{\L}{{\mathbb L}}

\newcommand{\x}{\mathbf{x}}
\newcommand{\y}{\mathbf{y}}
\newcommand{\e}{{\mathrm e}}

\newcommand{\n}{\mathbf n}

\newcommand{\hxi}{{\bm \xi}}

\begin{document}

\title{Accumulation time of diffusion in a 2D singularly perturbed domain}

\author{ \em
P. C. Bressloff, \\ Department of Mathematics, 
University of Utah \\155 South 1400 East, Salt Lake City, UT 84112}

\maketitle

\begin{abstract} 

A general problem of current interest is the analysis of diffusion problems in singularly perturbed domains, within which small subdomains are removed from the domain interior and boundary conditions imposed on the resulting holes. One major application is to intracellular diffusion, where the holes could represent organelles or biochemical substrates. In this paper we use a combination of matched asymptotic analysis and Green's function methods to calculate the so-called accumulation time for relaxation to steady state. The standard measure of the relaxation rate is in terms of the principal nonzero eigenvalue of the negative Laplacian. However, this global measure does not account for possible differences in the relaxation rate at different spatial locations, is independent of the initial conditions, and relies on the assumption that the eigenvalues have sufficiently large spectral gaps. As previously established for diffusion-based morphogen gradient formation, the accumulation time provides a better measure of the relaxation process. 

\end{abstract}

\section{Introduction}

In recent years there have been a large number of studies of diffusion problems in singularly perturbed domains, building on the seminal papers of Ward et. al. \cite{Ward93,Ward93a}. One general class of problem involves diffusion in a perforated domain, in which small subdomains are removed from the domain interior and boundary conditions imposed on the resulting holes \cite{Straube07,Bressloff08,Coombs09,Cheviakov11,Chevalier11,Ward15,Coombs15,Bressloff15,Lindsay16,Lindsay17,Grebenkov20,Bressloff21A,Bressloff21B}. Another important type of problem is diffusion in a domain with an exterior boundary that is reflecting almost everywhere, except for one or more small holes through which particles can escape \cite{Schuss07,Benichou08,Pillay10,Cheviakov10,Bressloff15a,Lindsay15,Holcman14a}. Much of the growth of interest in this area has been generated by models of diffusion in intracellular domains \cite{Holcman14b,Bressloff22}. For example, interior holes could represent subcellular structures such as organelles or biochemical substrates, whereas holes on the boundary of a domain could represent ion channels or nuclear pores. 
 Other applications include tracking the spread of chemical pollutants or heat from localized sources. 

Quantities of interest at the level of bulk diffusion include the steady-state solution (if it exists) and the approach to steady state as characterized by the leading non-zero eigenvalue $\lambda_1$ of the negative Laplacian. On the other hand, at the single-particle level, the diffusion equation represents the evolution of a probability density rather than a macroscopic particle concentration. One is now typically interested in the time for a particle to be captured by an interior trap (narrow capture) or to escape from a domain through a small hole in the boundary (narrow escape). Particular quantities of interest include splitting probabilities and conditional first passage times. In all of these cases, the resulting boundary value problems (BVPs) can be solved using a combination of matched asymptotic analysis and Green's function methods.

In this paper we are interested in the general problem of characterizing the relaxation to steady state of two-dimensional (2D) diffusion in a bounded domain with small interior holes. As we have already indicated, the standard approach developed by Ward {\em et. al.} \cite{Ward93,Ward93a} is to calculate the principal nonzero eigenvalue of the negative Laplacian. However, such a measure of the relaxation process has some limitations. First, it does not account for possible differences in the relaxation rate at different spatial locations. Second, it is independent of the initial conditions. Third, it relies on the assumption that the eigenvalues have sufficiently large spectral gaps. In light of this, we consider an alternative way of characterizing the approach to steady-state that is based on the the so-called accumulation time. The latter was originally developed within the context of diffusion-based morphogenesis \cite{Berez10,Berez11,Berez11a,Gordon11}, but has more recently been applied to intracellular protein gradient formation \cite{Bressloff19} and to diffusion processes with stochastic resetting \cite{Bressloff21C}. 

The structure of the paper is as follows. In section 2 we introduce the concept of an accumulation time by considering the simpler problem of one-dimensional (1D) diffusion-based concentration gradient formation. We also compare the accumulation time with the relaxation time based on an eigenfunction expansion. We then formulate the general problem of two-dimensional (2D) diffusion in a singularly-perturbed domain $\Omega$ and define the associated accumulation time $T(\x)$ in terms of the Laplace transform of the concentration, see section 3. The main part of the paper is the calculation of $T(\x)$. This is achieved by solving the diffusion equation in Laplace space using a combination of matched asymptotic analysis and Green's function methods (section 4). Since the 2D Green's function has a logarithmic singularity, $G(\x,\x')\sim -\ln|\x-\x'|$, it is natural to perform an asymptotic expansion in $\nu=-1/\ln \epsilon$ where $\epsilon$ is a small parameter that specifies the size of the holes relative to the size of the domain $\Omega$. In the case of a single hole, we compare the asymptotic expansion of the full accumulation time with a truncated version based on an eigenfunction expansion. We show that the difference between the two quantities is particularly significant in a neighborhood of the initial location of the diffusing particles. 
Finally, we illustrate the analysis by considering holes in the unit disc (section 5).
\section{Accumulation time of a 1D diffusion process}

In order to introduce the notion of an accumulation time for a diffusion process, consider the simple case of diffusion along the finite interval, $x\in [0,L]$, with a constant flux $J$ at the end $x=0$, a reflecting boundary at $x=L$, and a constant rate of degradation $k$ in the bulk of the domain. The concentration $u(x,t)$ evolves according to the equation 
\begin{equation}
\label{grad0}
\frac{\partial u}{\partial t}=-\L u\equiv D\frac{\partial^2u}{\partial x^2}-k u ,\ 0<x<L;\quad \left . -D\frac{\partial u}{\partial x}\right |_{x=0}=J, \quad \left . -D\frac{\partial u}{\partial x}\right |_{x=L}=0.
\end{equation}
Take the initial condition $u(x,0)=0$.
If $L \gg \xi\equiv \sqrt{{D}/{k}}$, then the boundary condition at $x=L$ can be ignored, and the steady-state solution is approximately given by a decaying exponential with space constant $\xi$:
\begin{equation}
\label{grad}
u^*(x)=\frac{J\xi}{D}\e^{-x/\xi}.
\end{equation}
One major application of the above type of model is to morphogen gradient formation during embryogenesis \cite{Wolpert69,Shvartsman12,Teimouri16}. In this particular case, $u(x,t)$ represents the extracellular morphogen concentration gradient along the body axis of a developing embryo, $k$ is an effective degradation or removal rate due to binding of morphogen to cell surface receptors, and the boundary flux at $x=0$ is generated by local protein synthesis. The spatially varying morphogen concentration drives a corresponding spatial variation in gene expression through some form of concentration thresholding mechanism. For example, in regions where the morphogen concentration exceeds a particular threshold, a specific gene is activated. Hence, a continuously varying morphogen concentration can be converted into a discrete spatial pattern of differentiated gene expression across a cell population. 
An important constraint on diffusion-based morphogenesis is that the concentration gradient across the length of the embryo should be established over appropriate developmental time scales. 

A standard method for estimating the time to approach steady state for diffusion in a bounded domain is to consider an eigenfunction expansion of the solution. That is, 
\begin{equation}
\label{eig}
u(x,t)-u^*(x)=\sum_{n\geq 0}c_n \phi_n(x)\e^{-\lambda_n t},
\end{equation}
where $0<\lambda_0< \lambda_1\ldots$ are the eigenvalues of the linear operator $\L=-D\partial^2_x+k$ with homogeneous boundary conditions, and the $\phi_n$ form a complete set of orthonormal eigenfunctions:
\[\L\phi_n(x) =\lambda_n\phi_n(x),\ x\in [0,L],\ \phi_n'(0)=\phi_n'(L),\quad \int_0^L\phi_m(x)\phi_n(x)dx=\delta_{n,m}.\]
The constant coefficients $c_n$ are determined by the initial condition. In the example given by equation (\ref{grad0}) one finds that
\begin{equation}
\phi_n(x)=A_n \cos(n\pi (L-x)/L),\quad \lambda_n = k+\frac{n^2\pi^2 D}{L^2}, \quad n\geq 0.
\end{equation}
If the positive eigenvalues are well separated, then the relaxation to steady state will be dominated by the term $c_0\phi_0(x)\e^{-\lambda_0t}$, and we can identify $1/\lambda_0$ as an effective relaxation time. For the given example, $\lambda_0=k$ and $\phi_0(x)=$ constant. However, characterizing the approach to steady state in terms of the smallest non-zero eigenvalue has some potential limitations. First, it does not account for possible differences in the relaxation rate at different spatial locations $x$. Second, all information regarding the initial condition is lost. Third, it relies on the assumption that the eigenvalues have sufficiently large spectral gaps, which may be difficult to establish in higher dimensions. In light of this, we will consider an alternative way of characterizing the approach to steady-state that is based on the the so-called accumulation time. The latter was originally developed within the context of diffusion-based morphogenesis \cite{Berez10,Berez11,Gordon11}, but has more recently been applied to intracellular protein gradient formation \cite{Bressloff19} and to diffusion processes with stochastic resetting \cite{Bressloff21C}. 

Let
\begin{equation}
\label{accu}
Z(x,t)=1-\frac{u(x,t)}{u^*(x)}
\end{equation}
be the fractional deviation of the concentration from steady state. Assuming that there is no overshooting, $1-Z(x,t)$ can be interpreted as the fraction of the steady-state concentration that has accumulated at $x$ by time $t$. It follows that $-\partial_t Z(x,t)dt$ is the fraction accumulated in the interval $[t,t+dt]$. The accumulation time $T(x)$ at position $x$ is then defined as \cite{Berez10,Berez11,Gordon11}:
\begin{equation}
\label{accu2}
T(x)=\int_0^{\infty} t\left (-\frac{\partial Z(x,t)}{\partial t}\right )dt=\int_0^{\infty} Z(x,t)dt.
\end{equation}
In terms of the eigenfunction expansion (\ref{eig}), 
\begin{equation}
T(x)=-\frac{1}{u^*(x)}\int_0^{\infty} \sum_{n\geq 0}c_n \phi_n(x)\e^{-\lambda_n t}dt=-\sum_{n=0}^{\infty}\frac{c_n \phi_n(x)}{\lambda_n u^*(x)},
\end{equation}
which is non-singular since $\lambda_n >0$ for all $n\geq 0$. If we simply kept the leading order term in the eigenvalue expansion, then $T(\x)\rightarrow T_0(\x)$ where
\begin{equation}
\label{Tapprox}
T_0(x)= -\frac{c_0\phi_0(x)}{\lambda_0 u^*(x)},
\end{equation}
Hence, approximating the solution using the leading-order term in an eigenfunction expansion does yield a position-dependent relaxation rate if one uses $T_0(x)$ rather than $\lambda_0^{-1}$ as the relaxation time. However, it can still be a poor approximation in regions where $\phi_0(x)\approx 0$ or if there is a small spectral gap. For the simple diffusion process given by equation (\ref{grad0}), $T(x)$ can be calculated explicitly in the limit $L\rightarrow \infty$ by considering the corresponding time-dependent solution \cite{Berez10}
\begin{eqnarray}
 u(x,t)=u^*(x)\left [1-\frac{1}{2}\mbox{erfc}\left (\frac{\sqrt{Dt}}{\xi}-\frac{x}{2\sqrt{Dt}}\right )-\frac{\e^{2x/\xi}}{2}\mbox{erfc}\left (\frac{\sqrt{Dt}}{\xi}+\frac{x}{2\sqrt{Dt}}\right )\right ],
\label{ct}
\end{eqnarray}
where $\mbox{erfc}(z)$ is the complementary error function. Substituting into equation (\ref{accu}) gives
\[Z(x,t)=\frac{1}{2}\mbox{erfc}\left (\frac{\sqrt{Dt}}{\xi}-\frac{x}{2\sqrt{Dt}}\right )+\frac{\e^{2x/\xi}}{2}\mbox{erfc}\left (\frac{\sqrt{Dt}}{\xi}+\frac{x}{2\sqrt{Dt}}\right ),\]
and, hence,
\begin{equation}
\label{morpht}
T(x)=\frac{1}{2k}\left (1+\sqrt{\frac{k}{D}}x\right ).
\end{equation}
On the other hand, in the limit $L\rightarrow \infty$, the spectral gap vanishes. Moreover, for sufficiently large $L$, equation (\ref{Tapprox}) yields
\begin{equation}
T_0(x)\approx \frac{A\sqrt{Dk}}{kJ}\e^{x\sqrt{k/D}},
\end{equation}
which is clearly a poor approximation of the exact accumulation time.

Finally, note that for the more complicated diffusion problems considered in this paper, it will be more convenient to calculate the accumulation time in Laplace space. Using the identity 
\[u^*(x)=\lim_{t\rightarrow \infty} u(x,t)=\lim_{s\rightarrow 0}s\widetilde{u}(x,s),\]
where $\u(x,s)=\int_0^{\infty}\e^{-st}u(x,t)dt$, and setting $\widetilde{F}(x,s)=s\widetilde{u}(x,s)$, the Laplace transform of equation (\ref{accu}) gives
\[s\widetilde{Z}(x,s)=1-\frac{\widetilde{F}(x,s)}{\widetilde{F}(x)},\quad \widetilde{F}(x)=\lim_{s\rightarrow 0}\widetilde{F}(x,s)=u^*(x)\]
and, hence
\begin{eqnarray}
 T(x)&=&\lim_{s\rightarrow 0} \widetilde{Z}(x,s) = \lim_{s\rightarrow 0}\frac{1}{s}\left [1-\frac{\widetilde{F}(x,s)}{\widetilde{F}(x)}\right ] =-\frac{1}{\widetilde{F}(x)}
\left .\frac{d}{ds}\widetilde{F}(x,s)\right |_{s=0}.
\label{acc}
\end{eqnarray}

\section{Diffusion in a 2D singularly perturbed domain}

\begin{figure}[b!]
\centering
\includegraphics[width=10cm]{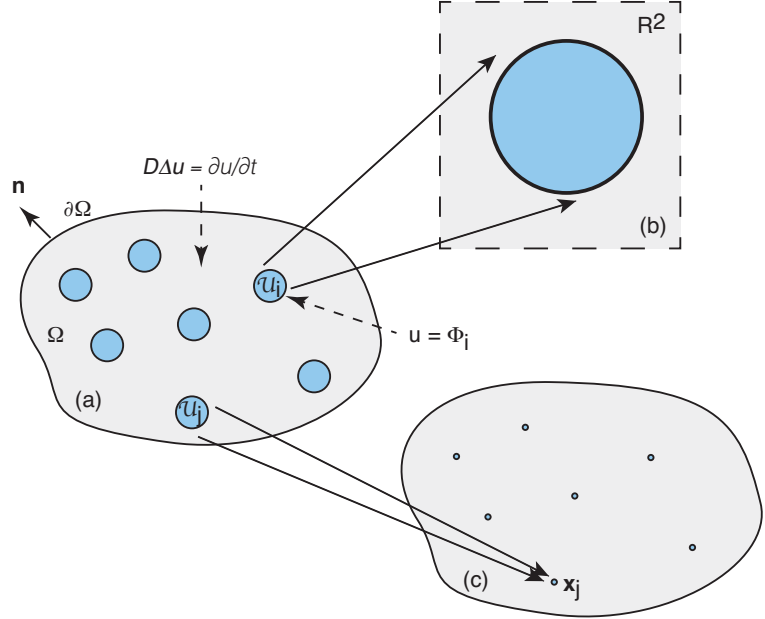} 
\caption{Diffusion in a 2D singularly perturbed domain. (a) Particles diffuse in a bounded domain $\Omega$ containing $N$ small interior holes or perforations denoted by $\calU_j$, $j=1,\ldots,N$. The exterior boundary $\partial \Omega$ is reflecting, whereas $u=\Phi_j$ on the $j$-th interior boundary $\partial \calU_i$. (b) Construction of the inner solution in terms of stretched coordinates $\y=\epsilon^{-1}(\x-{\x}_i)$, where ${\x}_i$ is the center of the $i$-th hole. The rescaled radius is $\rho_i=\ell_i$ and the region outside the hole is taken to be $\R^2$ rather than the bounded domain $\Omega$. (c) Construction of the outer solution. Each hole is shrunk to a single point. The outer solution can be expressed in terms of the corresponding modified Neumann Green's function and then matched with the inner solution around each hole.}
\label{fig1}
\end{figure}

Consider the diffusion equation in a bounded domain $\Omega\in \R^2$, that is perforated by a set of $N$ small holes denoted by $\calU_k\subset \Omega$, $k=1,\ldots,N$, see Fig. \ref{fig1}(a). The area of each hole is taken to be much smaller than $\Omega$, that is, $|\calU_j|\sim \epsilon^2 |\Omega|$ with $\calU_j\rightarrow \x_j\in \Omega$ uniformly as $\epsilon \rightarrow 0$, $j=1,\ldots,N$. In addition, the holes are assumed to be well separated with $|\x_i-\x_j|=O(1)$, $j\neq i$, and $\mbox{dist}(x_j,\partial \Omega)=O(1)$ for all $j=1,\ldots,N$. For simplicity, we take each hole to be a disc with $|\x-\x_j|=\epsilon \ell_j$. We impose a Neumann boundary condition on the external boundary $\partial \Omega$ and inhomogeneous Dirichlet boundary conditions on the interior boundaries $\partial \calU_j$. Let $u(\x,t)$ denote the concentration of freely diffusing particles for $\x\in \Omega\backslash \calU_a$, and $\calU_a\equiv\bigcup_{j=1}^N \calU_j$. Then
\begin{subequations} 
\label{RD1}
\begin{align}
	\frac{\partial u(\x,t)}{\partial t} &= D\nabla^2 u(\x,t),\ \x\in \Omega\backslash \calU_a,	\end{align}
together with the boundary conditions
\begin{equation}
\nabla u(\x,t)\cdot \n=0,\ \x \in \partial \Omega ; \quad u(\x,t)=\Phi_j, \ \x\in \partial \calU_j.
\end{equation}
Here $\n$ is the outward unit normal at a point on $\partial \Omega$. Finally, we impose the initial condition 
\begin{equation}
u(\x,0)=\Gamma_0\delta(\x-\x_0)
\end{equation}
 for some $\x_0 \in  \Omega\backslash \calU_a$, where $\Gamma_0$ is the initial number of molecules introduced into the domain.
 \end{subequations}
Equation (\ref{RD1}) represents one example of the general problem of diffusion in a singularly perturbed domain. Quantities of interest include the steady-state solution (if it exists) and the approach to steady state as characterized by the leading non-zero eigenvalue $\lambda_1$ of the negative Laplacian. The resulting boundary values problems (BVPs) can be solved using a combination of matched asymptotics and Green's function methods 
\cite{Ward93,Ward93a,Straube07,Bressloff08,Coombs09,Cheviakov11,Bressloff15,Grebenkov20}. A related class of BVPs arises when considering the capture of a single Brownian particle by small traps in the interior or boundary of the domain \cite{Schuss07,Chevalier11,Holcman14a,Coombs15,Ward15,Lindsay15,Lindsay16,Lindsay17,Bressloff21A,Bressloff21B}. In this case one is interested in moments of the conditional first passage time densities, for example. 

The main goal of this paper is to determine the accumulation time $T(\x)$ that characterizes the position-dependent relaxation to the steady state solution of equation (\ref{RD1}). This requires defining the higher-dimensional analogs of equations (\ref{accu}) and (\ref{acc}), namely,
\begin{equation}
\label{accuh}
Z(\x,t)=1-\frac{u(\x,t)}{u^*(\x)} 
\end{equation}
and
\begin{eqnarray}
T(\x)=\lim_{s\rightarrow 0} \widetilde{Z}(\x,s) = \lim_{s\rightarrow 0}\frac{1}{s}\left [1-\frac{\widetilde{F}(\x,s)}{\widetilde{F}(\x,0)}\right ] =-\frac{1}{\widetilde{F}(\x,0)}
\left .\frac{d}{ds}\widetilde{F}(\x,s)\right |_{s=0}
\label{acch}
\end{eqnarray}
with $\widetilde{F}(\x,s)=s\widetilde{u}(\x,s)$. By analogy with the example of morphogen gradient formation, we assume that $u(\x,t)<u^*(\x)$ for all $\x\in \Omega\backslash \calU_a$ and $ t >0$, that is, there is no overshooting. This can be ensured by taking
\begin{equation}
\frac{\Gamma_0}{|\Omega|} < \Phi_j \mbox{ for all }j=1,\ldots,N.
\label{grow}
\end{equation}
Equation (\ref{acch}) motivates solving the diffusion equation (\ref{RD1}) in Laplace space:
\begin{subequations}
\label{masterLT}
\begin{align}
 D\nabla^2\widetilde{u}-s\widetilde{u}  &=-\Gamma_0\delta(\x-\x_0), \quad \x \in \Omega\backslash \calU_a,\\
D\nabla \widetilde{u}(\x,s)\cdot \n&=0,\ \x \in \partial \Omega ,\quad  \widetilde{u}(\x,s)= \frac{\Phi_j}{s},\ \x\in \partial  \calU_j .
\end{align}
\end{subequations}
Note that one can eliminate the Dirac delta function on the right-hand side of equation (\ref{masterLT}a) by introducing the Green's function of the modified Helmholtz equation, 
\begin{align}
\label{GMH}
	D\nabla^2 G(\x,s|\x_0) -sG(\x,s|\x_0) &=-\delta(\x-\x_0) , \ \x\in \Omega,\quad
	 \nabla G(\x,s|x_0)\cdot \n  =0,\ \x \in \partial \Omega, 
	\end{align}
and taking
\begin{equation}
\label{uv}
\widetilde{u}(\x,s)=\Gamma_0G(\x,s|\x_0)+\v(\x,s),\ \x \in \Omega\backslash \calU_a,
\end{equation}
with
\begin{subequations} 
\label{masterLT2}
\begin{align}
	D\nabla^2 \v(\x,s ) -s\v(\x,s )&= 0 , \ \x\in \Omega\backslash \calU_a,\\
	 \nabla \v \cdot \n=0, \ \x\in \partial \Omega,\quad  \v&= \frac{\Phi_j}{s}-\Gamma_0G(\x,s|x_0),\ \x\in \partial\calU_i.
	\end{align}
\end{subequations}
We will derive an approximate solution of (\ref{masterLT}) by obtaining an inner or local solution valid in an $O(\epsilon)$ neighborhood of each hole, and then matching to an outer or global solution that is valid away from each neighborhood. 
The inner solution near the $j$-th hole is constructed by introducing the stretched local variable ${\mathbf y} = \varepsilon^{-1}(\x-\x_j)$ and setting
$U(\y,s)=\u(\x_j+\varepsilon \y,s)$, see Fig. \ref{fig1}(b). The resulting inner equation is 
\begin{align}
\label{inner}
& D\nabla^2_{\y}U = \epsilon^2 s U,\ |\y| > \ell_j,\quad U(\y,s)=\frac{\Phi_j}{s},\ |\y|=\ell_j.
\end{align}
The corresponding outer solution is constructed by shrinking each domain $\calU_j$ to a single point $\x_j$, see Fig. \ref{fig1}(c), and imposing a set of singularity conditions for $\x\rightarrow \x_j$, $j=1,\ldots,N$. The latter are determined by matching with the far field behavior of the $N$ inner solutions.

\section{Matched asymptotic analysis of the accumulation time}

A well known feature of diffusion in 2D singularly perturbed domains is that the matching of inner and outer solutions leads to an asymptotic expansion in powers of the small parameter $\nu=-1/\ln \epsilon$ rather than $\epsilon$ itself. This reflects the singular nature of the corresponding modified Helmholtz Green's function. That is, as $|\x-\x_0|\rightarrow 0$, we have
\begin{equation}
G(\x,s|\x_0)\rightarrow -\frac{1}{2\pi D}\ln|\x-\x_0|.
\end{equation}
Therefore, one typically considers $O(1)$ accuracy with respect to an $\epsilon$ expansion\footnote{Certain care has to be taken in considering asymptotic expansions with respect to $\nu$, since $\nu \rightarrow 0$ more slowly than $\epsilon\rightarrow 0$. Therefore, if one is interested in obtaining numerically accurate solutions at $O(1)$ in $\epsilon$,  then it is often necessary to sum over the logarithmic terms non-perturbatively along the lines of Ward and Keller \cite{Ward93}. This is equivalent to calculating the asymptotic solution for all terms of $O(\nu^k)$ for any $k$. In this paper, we will use the non-perturbative solution to generate the asymptotic expansion in $\nu$.}. At this level of approximation we can neglect 
 all terms of $O(\epsilon)$, so that the inner equation (\ref{inner}) becomes
 \begin{align}
\label{inner2D}
& D\nabla^2_{\y}U =0,\ |\y| >1,\quad U(\y,s)= \frac{\Phi_j}{s},\ |\y|=1.
\end{align}
We have also set $\ell_j=1$ for all $j$.
Using polar coordinates, the inner solution can be written as
\begin{align}
    &U = \Phi_j+\nu D^{-1}A_j(\nu, s) \log\rho  ,\quad 1\leq \rho < \infty.
\end{align}
The corresponding outer solution is given by equation (\ref{uv}) with 
\begin{subequations} 
\label{outer}
\begin{align}
	D\nabla^2 \v(\x,s ) -s\v(\x,s)&=0 , \ \x\in \Omega \backslash \{\x_1,\ldots,\x_N\},\quad
	 \nabla \v\cdot \n=0, \ \x\in \partial \Omega,	\end{align}
together with the $N$ singularity conditions
	 \begin{equation}
	 \v(\x,s)\sim V_j(s) +\nu D^{-1}A_j(\nu,s)\ln|\x-\x_j|/\epsilon \quad \mbox{as} \ \x\rightarrow \x_j,\quad j=1,\ldots,N,
	 \label{sing}
	 \end{equation}
\end{subequations}
with
\begin{equation}
\label{Vj}
V_j(s)=\frac{\Phi_j}{s}-\Gamma_0G(\x_j,s|\x_0).
\end{equation}
Therefore,	 the outer solution can be written as
\begin{equation}
\label{out}
\v(\x,s )=-2\pi \nu  \sum_{j=1}^N A_j(\nu,s)G(\x,s|\x_j),
\end{equation}
where $G$ is the Green's function defined in equation (\ref{GMH}). It is convenient to isolate the logarithmic singularity of $G$ by defining the regular part $R$:
\begin{equation}
G(\x,s|\x_0)=-\frac{1}{2\pi D}\ln|\x-\x_0|+R(\x,s|\x_0),
\end{equation}

There remain $N$ unknown coefficients $A_j(\nu,s)$, which require $N$ additional constraints. The latter are obtained by matching the near-field behavior of the outer solution (\ref{out}) with the singularity condition (\ref{sing}) in a neighborhood of $\calU_j$ for $j=1,\ldots,N$:
\begin{equation}
\label{match1}
2\pi \nu \sum_{i\neq j}  A_{i}(\nu,s)G(\x_i,s|\x_j)+D^{-1}A_{j}(\nu,s)+2\pi \nu  A_{j}(\nu,s)R(\x_j,s|\x_j)
= -V_j .
\end{equation}
The matrix equation (\ref{match1}) has the formal solution
\begin{equation}
\label{matrix}
A_{i}(\nu,s)=-D\sum_{j=1}^N[\delta_{i,j}+2\pi \nu D {\mathcal G}_{ij}(s)]^{-1}V_j,\quad i=1,\ldots,N,
\end{equation}
where ${\mathcal G}_{ij}(s)= G(\x_i,s|\x_j)$ for $ i\neq j$ and ${\mathcal G}_{ii}(s)= R(\x_i,s|\x_i)$.
Equations (\ref{out}) and (\ref{matrix}) yield a non-perturbative solution that sums over all logarithmic terms involving factors of $\nu$, along analogous lines to \cite{Ward93}; the solution is $O(1)$ with respect to a corresponding $\epsilon$ expansion. (One could also include $O(\epsilon)$ and higher-order terms as illustrated in Ref. \cite{Lindsay15}; here we will focus on the $O(1)$ expansion.)

\subsection{Steady-state solution}

Multiplying equation (\ref{out}) by $s$ and then taking the limit $s\rightarrow 0$ yields the steady-state solution
\begin{align}
\label{split2}
v^*(\x)&=-2\pi \nu \lim_{s\rightarrow 0}s\sum_{j=1}^N A_j(\nu,s)G(\x,s|\x_j).
\end{align}
In order to calculate the above limit, we use the result that
\begin{equation}
\label{Gs}
G(\x,s|\x_0)=\frac{1}{s|\Omega|}+G_0(\x,\x_0)+sG_1(\x,\x_0)+O(s^2),
\end{equation}
where $G_0$ is the generalized Neumann Green's function of Laplace's equation:
\begin{subequations}
\label{GM}
\begin{align}
	D\nabla^2 G_0(\x,\x_0)  &=\frac{1}{|\Omega|}-\delta(\x-\x_0) , \ \x\in \Omega,\\
	 \nabla G_0(\x,\x_0)\cdot \n  &=0,\ \x \in \partial \Omega,\quad \int_{\calU}G_0(\x,\x_0)d\x=0, \\
	 G_0(\x,\x_0)&=-\frac{1}{2\pi D}\ln|\x-\x_0|+R_0(\x,\x_0).
	\end{align}
	\end{subequations}
It follows that the coefficient $V_j$ has the small-$s$ expansion
\begin{equation}
V_j=\frac{\Phi_j-\Gamma_0/|\Omega|}{s}-\Gamma_0G_0(\x_j,\x_0)+O(s).
\end{equation}
and we can introduce an analogous small-$s$ expansion for the coefficient $A_k(\nu,s)$:
\begin{equation}
\label{anu}
A_k(\nu,s)=\frac{{A}_k(\nu)}{s}+\theta_k(\nu)+s\chi_k(\nu)+O(s^2)
\end{equation}
Equations (\ref{split2}), (\ref{Gs}) and (\ref{anu}) imply that
\begin{align}
v^*(\x)&=-2\pi \nu \lim_{s\rightarrow 0}s\sum_{k=1}^N \left [ \frac{{A}_k(\nu)}{s}+\theta_k(\nu)+O(s)\right ] \left [\frac{1}{s|\Omega|}+G_0(\x,\x_k)+O(s)\right ]\nonumber \\
&=-2\pi \nu\sum_{k=1}^N\left (\frac{\theta_k(\nu)}{|\Omega|}+{A}_k(\nu)G_0(\x,\x_k)\right ).
\end{align}

The unknown coefficients ${A}_k(\nu)$ and $\theta_k(\nu)$ are determined by substituting the various small-$s$ expansions into equation (\ref{match1}): 
\begin{align}
&2\pi \nu \sum_{i=1}^N \left [\frac{{A}_i(\nu)}{s}+\theta_i(\nu)+s\chi_i\nu)+O(s^2)\right ]\left (\frac{1}{s|\Omega|}+{\mathcal G}_{ij}^{(0)}+s{\mathcal G}_{ij}^{(1)}+O(s^2)\right )
\nonumber\\
&\qquad +D^{-1}\left [\frac{{A}_j(\nu)}{s}+\theta_j(\nu)+O(s)\right ]=\frac{\Gamma_0/|\Omega|-\Phi_j}{s}+\Gamma_0G_0(\x_j,\x_0)+O(s),
\label{match2}
\end{align}
where ${\mathcal G}_{ij}^{(n)}= G_{n}(\x_i,\x_j)$ for $ i\neq j$ and ${\mathcal G}_{ii}^{(n)}= R_{n}(\x_i,\x_i)$, $n=0,1$.
Collecting the $O(s^{-2})$ and $O(s^{-1})$ terms yields the pair of conditions
\begin{subequations}
\label{cons}
\begin{align}
&\sum_{i=1}^N{A}_i(\nu)=0,\\
&2\pi \nu A(\nu)\circ {\mathcal G}^{(0)}_j+\frac{2\pi \nu }{|\Omega|}\sum_{i=1}^N{\theta_i}(\nu)+\frac{{A}_j(\nu)}{D}=-\Phi_j+\frac{\Gamma_0}{|\Omega|},\quad j=1,\ldots,N.
\end{align}
\end{subequations}
We have introduced the notation
\begin{equation}
[A(\nu)\circ {\mathcal G}^{(n)}]_j=\sum_{i=1}^N{A}_i(\nu){\mathcal G}_{ij}^{(n)}, \quad n\geq 0,
\end{equation}
and similarly for any other $N$-vector.
Finally, imposing the condition
\begin{equation}
\label{Asum}
 \sum_{i=1}^N{\theta}_i(\nu)=\frac{\Gamma_0+\Delta \Gamma}{2\pi \nu}
\end{equation}
for some unknown constant $\Delta \Gamma$, and using equation (\ref{uv}), yields the steady-state solution
\begin{equation}
u^*(\x)=-\frac{\Delta \Gamma}{|\Omega|}-2\pi \nu\sum_{k=1}^N {A}_k(\nu)G_0(\x,\x_k)
\label{ustar}
\end{equation}
with the $N+1$ unknowns $\Delta \Gamma$ and ${A}_k(\nu),k=1,\ldots,N$ satisfying the $N+1$ equations
\begin{subequations}
\label{Achi}
\begin{align}
&\sum_{i=1}^N{A}_i(\nu)=0,\\
&2\pi \nu [ A(\nu)\circ {\mathcal G}^{(0)}]_j+\frac{{A}_j(\nu)}{D}=-\Phi_j-\frac{\Delta \Gamma}{|\Omega|}.
\end{align}
\end{subequations}
Summing both sides of equation (\ref{Achi}b) with respect to $j$ and using (\ref{Achi}a) gives (after dividing through by $N$)
\begin{align}
\label{chi}
-\frac{\Delta \Gamma}{|\Omega|}&=\overline{\Phi}+2\pi \nu \overline{A(\nu)\circ {\mathcal G}^{(0)}}
\end{align}
with
\begin{equation}
\overline{\Phi}=\frac{1}{N}\sum_{j=1}^N\Phi_j,\quad \overline{A(\nu)\circ {\mathcal G}^{(n)}}=\frac{1 }{N}\sum_{i,j=1}^N {A}_i(\nu){\mathcal G}_{ij}^{(n)}.
\end{equation}
We can thus eliminate $\Delta \Gamma$ from equations (\ref{ustar}) and (\ref{Achi}b) such that 
\begin{equation}
u^*(\x)=\overline{\Phi}-2\pi \nu\sum_{k=1}^N {A}_k(\nu)G_0(\x,\x_k)
+2\pi \nu \overline{A(\nu)\circ {\mathcal G}^{(0)}}
\label{ustar2}
\end{equation}
with the $N-1$ independent coefficients $A_i$ determined from the matrix equation
\begin{align}
\label{match3}
2\pi \nu [A(\nu)\circ {\mathcal G}^{(0)}]_j-2\pi \nu \overline{A(\nu)\circ {\mathcal G}^{(0)}}+\frac{{A}_j(\nu)}{D}=\overline{\Phi}-\Phi_j.
\end{align}
Finally, note that we could have obtained the steady-state solution (\ref{match3}) more directly by applying matched asymptotics to the steady-state version of equation (\ref{RD1}).  However, working in Laplace space also allows us to calculate the accumulation time.

\subsection{Accumulation time}

In order to calculate the accumulation time $T(\x)$ given by equation (\ref{acch}), we 
need to determine the first derivative 
${\mathcal F}(\x)\equiv \left .{d\widetilde{F}(\x,s)}/{ds}\right |_{s=0}$
with $\widetilde{F}(\x,s)=s\widetilde{u}(\x,s)$ and $\u(\x,s)$ given by equation (\ref{uv}) and the outer solution (\ref{out}):
\begin{equation}
\u(\x,s )\sim\Gamma_0G(\x,s|\x_0)-2\pi \nu  \sum_{j=1}^N A_j(\nu,s)G(\x,s|\x_j).
\end{equation}
Using the small-$s$ expansions (\ref{Gs}) and (\ref{anu}), we have
\begin{align}
{\mathcal F}(\x)&\sim\Gamma_0G_0(\x,\x_0)-2\pi \nu \lim_{s\rightarrow 0} \frac{d}{ds} s\sum_{j=1}^N\left ( \frac{{A}_j(\nu)}{s}+\theta_j(\nu)+s\chi_j(\nu)+O(s^2)\right )\nonumber \\
&\hspace{4cm} \times \left (\frac{1}{s|\Omega|}+G_0(\x,\x_j)+sG_1(\x,\x_j)+O(s^2)\right )\nonumber \\
&\sim\Gamma_0G_0(\x,\x_0)-2\pi \nu \sum_{j=1}^N\left ( {A}_j(\nu) G_1(\x,\x_j) +
{\theta}_j(\nu) G_0(\x,\x_j)+\frac{\chi_j(\nu)}{|\Omega|}\right ).
\label{calF}
\end{align}
Since ${\mathcal F}(\x)$ involves the higher-order term $\sum_{k=1}^N\chi_k$, it is necessary to consider the $O(1)$ contributions to equation (\ref{match2}):
\begin{align}
2\pi\nu[\theta(\nu)\circ {\mathcal G}^{(0)}]_j+2\pi\nu  [A(\nu)\circ {\mathcal G}^{(1)}]_j+\frac{2\pi \nu}{|\Omega|}\sum_{i=1}^N{\chi}_i(\nu) +\frac{{\theta}_j(\nu)}{D}=\Gamma_0G_0(\x_j,\x_0)
\end{align}
for $ j=1,\ldots,N$.
Summing both sides with respect to $j$ implies that
\begin{align}
\frac{2\pi \nu}{|\Omega|}\sum_{i=1}^N{\chi}_i(\nu) =\frac{\Gamma_0}{N}\sum_{j=1}^NG_0(\x_j,\x_0)-2\pi \nu \overline{\theta(\nu)\circ {\mathcal G}^{(0)}}-2\pi \nu \overline{A(\nu)\circ {\mathcal G}^{(1)}}-\frac{1}{D}\frac{\Gamma_0+\Delta \Gamma}{2\pi \nu N}.
\label{bobo}
\end{align}
We have used the identity (\ref{Asum}). Combining equations (\ref{chi}) and (\ref{bobo}) and substituting into (\ref{calF}) gives 
\begin{align}
{\mathcal F}(\x)&\sim\frac{1}{D}\frac{\Gamma_0-|\Omega|\overline{\Phi}}{2\pi \nu N}+\Gamma_0G_0(\x,\x_0)-\frac{\Gamma_0}{N}\sum_{j=1}^NG_0(\x_j,\x_0)-\frac{|\Omega| }{ND}\overline{A(\nu)\circ {\mathcal G}^{(0)}}\ \\
&\quad +2\pi \nu \overline{\theta(\nu)\circ {\mathcal G}^{(0)}}+2\pi \nu \overline{A(\nu)\circ {\mathcal G}^{(1)}}-2\pi \nu \sum_{j=1}^N\left ( {A}_j(\nu) G_1(\x,\x_j) +
{\theta}_j(\nu) G_0(\x,\x_j)\right ).\nonumber \end{align}

Finally, using equation (\ref{acch}) and the steady-state solution (\ref{ustar2}), we obtain the following non-perturbative expression for the accumulation time:
\begin{align}
\label{TT0}
T(\x)&=-\frac{{\mathcal F}(\x)}{u^*(\x)}=-\frac{{\mathcal F}(\x)}{\overline{\Phi}-2\pi \nu\sum_{k=1}^N {A}_k(\nu)G_0(\x,\x_k)
+2\pi \nu \overline{A(\nu)\circ {\mathcal G}^{(0)}}}.
\end{align}
We now note that
\begin{equation}
A_j(\nu)=D(\overline{\Phi}-\Phi_j)+O(\nu),\quad 2\pi \nu \theta_j(\nu)=\frac{\Gamma_0-|\Omega|\overline{\Phi}}{N}+O(\nu).
\end{equation}
Hence,
\begin{align}
{\mathcal F}(\x)\sim  \frac{\Gamma_0-|\Omega|\overline{\Phi}}{2\pi \nu ND}+{\mathcal F}_0(\x)+O(\nu)
\end{align}
with
\begin{align}
{\mathcal F}_0(\x)&=\Gamma_0G_0(\x,\x_0)-\frac{\Gamma_0}{N}\sum_{j=1}^NG_0(\x_j,\x_0)-\frac{|\Omega| }{N^2}\sum_{i,j=1}^N(\overline{\Phi}-\Phi_i){\mathcal G}_{ij}^{(0)}\nonumber\\
&\quad +\frac{\Gamma_0-|\Omega|\overline{\Phi}}{N^2}\sum_{i,j=1}^N\calG_{ij}^{(0)}-\frac{\Gamma_0-|\Omega|\overline{\Phi}}{N}\sum_{j=1}^NG_0(\x,\x_j) .
\label{F02D}
\end{align}
Similarly, the steady-state density given by equation (\ref{ustar2}) has the $\nu$-expansion
\begin{equation}
u^*(\x)=\overline{\Phi}-2\pi \nu\left [\sum_{k=1}^N (\overline{\Phi}-\Phi_k)G_0(\x,\x_k)
-\frac{1}{N}\sum_{j,k=1}^N (\overline{\Phi}-\Phi_j) {\mathcal G}_{jk}^{(0)}\right ]+O(\nu^2).
\label{ustar3}
\end{equation}
We thus obtain the following $\nu$-expansion for the accumulation time to $O(\nu)$:
\begin{align}
T(\x)&\sim\frac{1}{D}\frac{|\Omega|\overline{\Phi}-\Gamma_0}{2\pi \nu N\overline{\Phi}}-\frac{{\mathcal F}_0(\x)}{\overline{\Phi}}+\frac{|\Omega|\overline{\Phi}-\Gamma_0}{  N\overline{\Phi}^2}\left [\sum_{k=1}^N (\overline{\Phi}-\Phi_k)G_0(\x,\x_k)
-\frac{1}{N}\sum_{j,k=1}^N (\overline{\Phi}-\Phi_j) {\mathcal G}_{jk}^{(0)}\right ].
\label{Tres}
\end{align}

A number of results follow from the above asymptotic analysis.
\medskip

\noindent (i) The leading order contribution to $T(\x)$ is the constant $\mu_0/\nu$ with
$\mu_0=(|\Omega|\overline{\Phi}-\Gamma_0)/{2\pi  ND\overline{\Phi}}$.
Note that $\mu_0>0$ due to the condition (\ref{grow}), which ensures that the accumulation time is positive.
Moreover, $T(\x)\rightarrow \infty$ as $\nu \rightarrow 0$. This singular behavior as the size of the holes shrinks to zero is related to the fact that
$\lim_{\nu \rightarrow 0}u^*(\x)=\overline{\Phi}$,
whereas the steady-state in the absence of any holes is $\Gamma_0/|\Omega|$. In other words, the limits $\nu\rightarrow 0$ and $t\rightarrow \infty$ do not commute.
\medskip

\noindent (ii) In the case of $N$ identical interior boundary conditions, $\Phi_j=\overline{\Phi}$ for all $j$, the expression for the accumulation time simplifies according to
\begin{align}
\label{one}
T(\x)&=\frac{1}{D}\frac{|\Omega|\overline{\Phi}-\Gamma_0}{2\pi \nu N \overline{\Phi}}-\frac{\Gamma_0}{\overline{\Phi}}\left [ G_0(\x,\x_0)-\frac{1}{N}\sum_{j=1}^NG_0(\x_j,\x_0)\right ]\\
&\quad -\frac{|\Omega|\overline{\Phi}-\Gamma_0}{ N \overline{\Phi}}\bigg \{  \sum_{j=1}^NG_0(\x,\x_j)-\frac{1}{N}\sum_{i,j=1}^N\calG_{ij}^{(0)}\bigg \}+O(\nu).\nonumber 
\end{align}

\noindent (iii) In general, the accumulation time depends on the hole positions $\x_j$ at $O(1)$ with respect to a power series expansion in $\nu$. On the other hand, this position dependence occurs at $O(\nu)$ in the $\nu$-expansion of equation (\ref{ustar2}) for the steady state $u^*(\x)$. That is, the accumulation time is a more sensitive function of the spatial distribution of the holes. 
\medskip

\noindent (iv) The accumulation time is strongly dependent on the initial condition $u(\x,0)=\Gamma_0\delta(\x-\x_0)$ via its dependence on $\Gamma_0$ and $\x_0$. Note that the singularity in $T(\x)$ as $\x\rightarrow \x_0$ is a consequence of the Dirac delta function. It is easily removed by taking the initial concentration to be a strongly localized Gaussian, for example.

\subsection{Eigenfunction expansion}

Characterizing the relaxation to steady state in terms of the $\x$-dependent accumulation time $T(\x)$ is significantly different from the standard method based on an eigenvalue expansion \cite{Ward93,Coombs09}. Consider the set of eigenpairs of the negative Laplacian in the given singularly perturbed domain, which are denoted by $(\lambda_n,\phi_n(\x))$ for $n\geq 0$ with $0<\lambda_0<\lambda_1 <\lambda_2\ldots$ and 
\[\int_{\Omega\backslash \calU_a}\phi_n(\x)\phi_m(\x)d\x=\delta_{n,m}.\]
 Then
\begin{equation}
\label{eig2}
u(\x,t)-u^*(\x)=\sum_{n\geq 0}c_n \phi_n(\x)\e^{-\lambda_n t}\approx c_0\phi_0(\x)\e^{-\lambda_0 t},
\end{equation}
where $\lambda_0$ is the smallest nonzero eigenvalue. Since $\lambda_0=O(\nu)$, it can be calculated by solving the singularly perturbed BVP \cite{Ward93,Coombs09} 
 \begin{subequations}
 \label{BVPeig}
\begin{align}
&D\nabla^2 \phi+\lambda\phi=0,\ \x\in \Omega\backslash \calU_a,\quad \nabla\phi \cdot \n =0,\ \x \in \partial \Omega, \ \int_{\Omega\backslash \calU_a}\phi^2(\x)d\x=1,\\
&\phi=0,\ \x\in \partial \calU_j,\quad j=1,\ldots,N.
\end{align}
\end{subequations}
(We have dropped the subscripts on the leading eigenvalue and eigenfunction.)
In particular, the outer solution for the principal eigenfunction satisfies the equation
 \begin{subequations}
\begin{align}
&D\nabla^2 \phi+\lambda\phi=0,\ \x\in \Omega\backslash \{\x_1,\ldots,\x_N\},\quad \nabla\phi\cdot \n =0,\ \x \in \partial \Omega, \ \int_{\Omega}\phi^2(\x)d\x=1,\\
&\phi \sim \nu D^{-1}B_j(\nu,\lambda)\ln|\x-\x_j|/\epsilon \quad \mbox{as} \ \x\rightarrow \x_j,\quad j=1,\ldots,N.
\end{align}
\end{subequations}
It follows that
\begin{equation}
\label{outphi}
\phi(\x,s )=-2\pi \nu  \sum_{j=1}^N B_j(\nu,\lambda)G(\x,-\lambda|\x_j),
\end{equation}
where $G$ is the 2D modified Helmholtz Green's function defined in equation (\ref{GMH}). The $N$ unknown coefficients $B_j(\nu,\lambda)$ are obtained by matching the far-field behavior of the corresponding inner solution with the near-field behavior of the outer solution in a neighborhood of $\calU_j$ for $j=1,\ldots,N$:
\begin{equation}
\label{Amatch1}
2\pi \nu \sum_{i=1}^N  B_{i}(\nu,\lambda){\mathcal G}_{ij}(-\lambda)+D^{-1}B_{j}(\nu,\lambda)
=0.
\end{equation}
we thus obtain a transcendental equation for $\lambda_1$ of the form
\begin{equation}
\label{det}
\mbox{Det}\left (2\pi \nu D \bm{\mathcal G}(-\lambda)+{\bf I}\right )=0,
\end{equation}
where ${\bf I}$ is the $N\times N$ unit matrix.

In the case of a single hole ($N=1$) equation (\ref{det}) simplifies to the condition
\begin{equation}
\label{bbtrans}
R(\x_1,-\lambda|\x_1)=-\frac{1}{2\pi \nu D}.
\end{equation}
Since $\lambda=\nu \lambda_1+\nu^2 \lambda_2+\ldots$ we can use the Green's function expansion (\ref{Gs}):
\begin{equation*}
G(\x_i,-\lambda|\x_j)=-\frac{1}{\lambda|\Omega|}+G_0(\x_i,\x_j)-\lambda G_1(\x_i,\x_j)+O(\lambda^2).
\end{equation*}
Substituting the corresponding expansion for the regular part of the Green's function into equation (\ref{Amatch1}) and performing an asymptotic expansion in $\nu$ we find that \begin{equation}
\lambda = \frac{2\pi \nu D }{|\Omega|} \left (1- {2\pi \nu D} R_0(\x_1,\x_1) \right ) +O(\nu^3),
\end{equation}
so that the relaxation time is
\begin{equation}
\tau\equiv \frac{1}{\lambda}= \frac{|\Omega|}{2\pi \nu D }  \left (1+2\pi \nu D R_0(\x_1,\x_1)\right ) +O(\nu).
\end{equation}
Comparison with equation (\ref{one}) for $N=1$ shows that $\tau$ is a global rather than a local measure of relaxation to the steady state, and is independent of the initial density and the boundary values $\Phi_j$. As expected, $\tau\rightarrow \infty$ as $\nu\rightarrow 0$. Now consider the approximation of the accumulation time given by the 2D analog of equation (\ref{Tapprox}).
First, the principal eigenfunction is
\begin{equation}
\phi(\x)= A G(\x,-\lambda|\x_1) \sim A \left (-\frac{\tau}{ |\Omega|}+G_0(\x,\x_1)\right ),
\end{equation}
with the constant $A$ determined by the normalization condition $\int_{\Omega}\phi^2(\x)d\x =1$. That is,
\begin{align}
1=A^2 \int_{\Omega}\left (\frac{\tau^2}{|\Omega|^2}-\frac{2\tau}{|\Omega|}G_0(\x,\x_1)+O(1)\right )d\x,
\end{align}
which implies
\begin{equation}
A\sim \frac{\sqrt{|\Omega|}}{\tau}.
\end{equation}
Second, the coefficient $c_0$ is obtained by setting $t=0$ in equation  (\ref{eig2}):
\begin{equation}
\Gamma_0\delta(\x-\x_0) -\Phi_1 =\sum_{n=0}^{\infty}c_n\phi_n(\x).
\end{equation}
Note that $u^*(\x)=\Phi_1$ for a single hole. 
Multiplying both sides by $\phi_0(\x)=\phi(\x)$, integrating with respect to $\x$ and imposing orthonormality of the eigenfunctions yields
\begin{equation}
c_0=\Gamma_0\phi(\x_0)-\Phi_1\int_{\Omega}\phi(\x)d\x.
\end{equation}
Substituting the solution for $\phi(\x)$ and using the normalization
condition $\int_{\Omega}G(\x,s|\x_0)d\x=1/s$, we have
\begin{equation}
c_0=A \bigg (\Gamma_0G(\x_0,-\lambda|\x_1) +\Phi_1 \tau \bigg )\sim A\left \{\left (\Phi_1-\frac{\Gamma_0}{|\Omega|}\right )\tau+\Gamma_0 G_0(\x_0,\x_1)\right \}.
\end{equation}
Finally, combining our various results shows that
\begin{align}
\label{Tapprox2}
T_0(x)&=-\frac{c_0\phi(x)}{\lambda \Phi_1}\sim \frac{|\Omega|}{\Phi_1 \tau}\left \{\left (\Phi_1-\frac{\Gamma_0}{|\Omega|}\right )\tau+\Gamma_0 G_0(\x_0,\x_1)\right \}\left (\frac{\tau}{ |\Omega|}-G_0(\x,\x_1)\right )\nonumber \\
&\sim \frac{|\Omega| {\Phi_1}-\Gamma_0}{2\pi \nu D {\Phi_1}}+\frac{\Gamma_0}{\Phi_1} G_0(\x_1,\x_0)+\frac{|\Omega| \Phi_1-\Gamma_0}{\Phi_1}[R_0(\x_1,\x_1)-G_0(\x,\x_1) ] +O(\nu).
\end{align}
On the other hand, setting $N=1$ in equation (\ref{one}) shows that the asymptotic expansion of the full accumulation time is
\begin{align}
T(\x)&\sim T_0(\x)-\frac{\Gamma_0}{\Phi_1}G_0(\x,\x_0)+O(\nu).
\label{Tres1}
\end{align}

\section{Examples} 

\begin{figure}[b!]
\centering
\includegraphics[width=12cm]{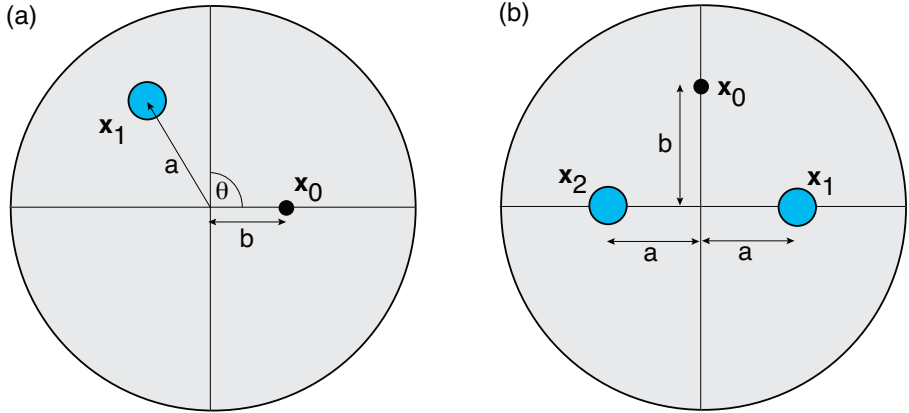} 
\caption{Example target configurations. (a) A single hole of radius $\epsilon$ is placed at position $\x_1=a(\cos \theta,\sin \theta)$ in the unit disc with $a<1$.The initial condition is taken to be a localized concentration on the $x$-axis, $u(\x,0)=\Gamma_0 \delta(\x-\x_0)$ with $\x_0=(0,b)$ and $0<n<1$ $\x_0=(b,0)$.  (b) A pair of identical small holes are placed at the points $\x_{1,2} =(\pm a,0)$ along the $x$-axis of the unit disc with $0<a <1$. The boundary conditions are $\partial_nu(x,t)=0$ for $|\x|=1$ and $u(\x,t)=1$ for $|\x-{\bf x}_{1,2}|=\epsilon$. The initial position $\x_0$ is now taken to be on the $y$-axis with $\x_0=(0,b)$, $0 < b < 1$.}
\label{fig2}
\end{figure}

\begin{figure}[b!]
\centering
\includegraphics[width=13cm]{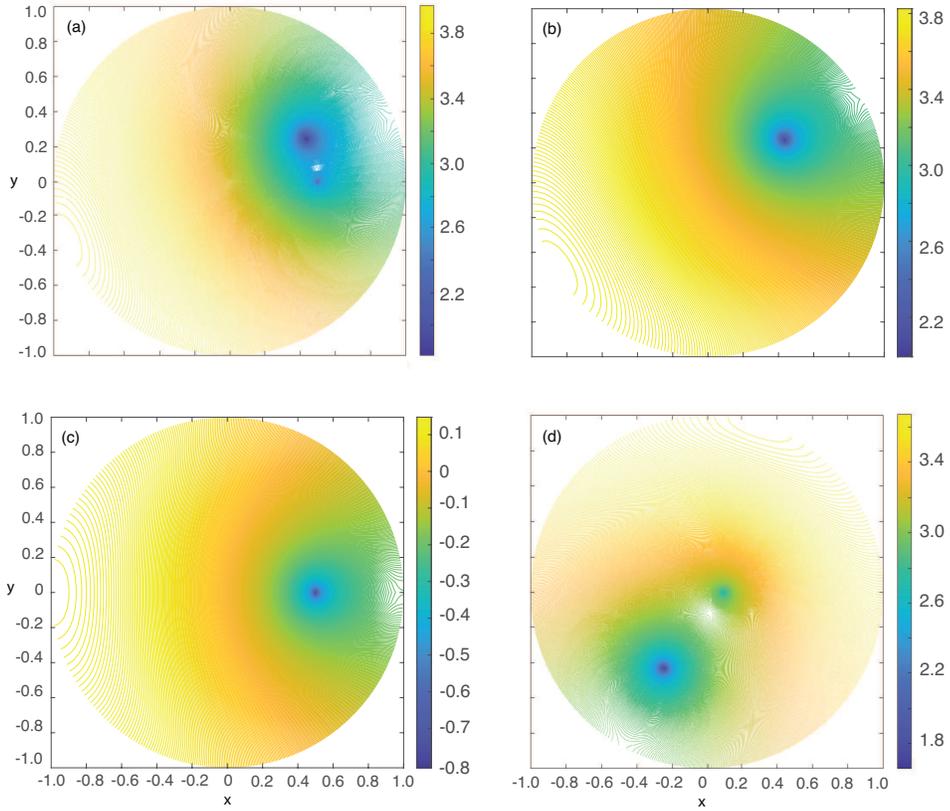} 
\caption{Accumulation time in the unit disc for the single-hole configuration shown in Fig. \ref{fig2}(a) with $\x_1=a(\cos \theta_1,\sin\theta_1)$ and $\x_0=(b,0)$. (a) Contour plot of the full accumulation time $T(\x)$ for $a=b=0.5$ and $\theta_1=\pi/6$. (b,c) Corresponding plots of the truncated accumulation time $T_0(\x)$ and the difference $T(\x)-T_0(\x)$. (d) Plot of $T(\x)$ as a function of $r$ for $\x=r(\cos \theta,\sin \theta)$, $a=0.5$, $b=0.1$ and $\theta_1=4\pi/3$. Other parameter values are $\Gamma_0=1$, $\nu=0.1$ and $D=1$.}
\label{fig3}
\end{figure}

\subsection{Single target in the unit disc}

\begin{figure}[t!]
\centering
\includegraphics[width=13cm]{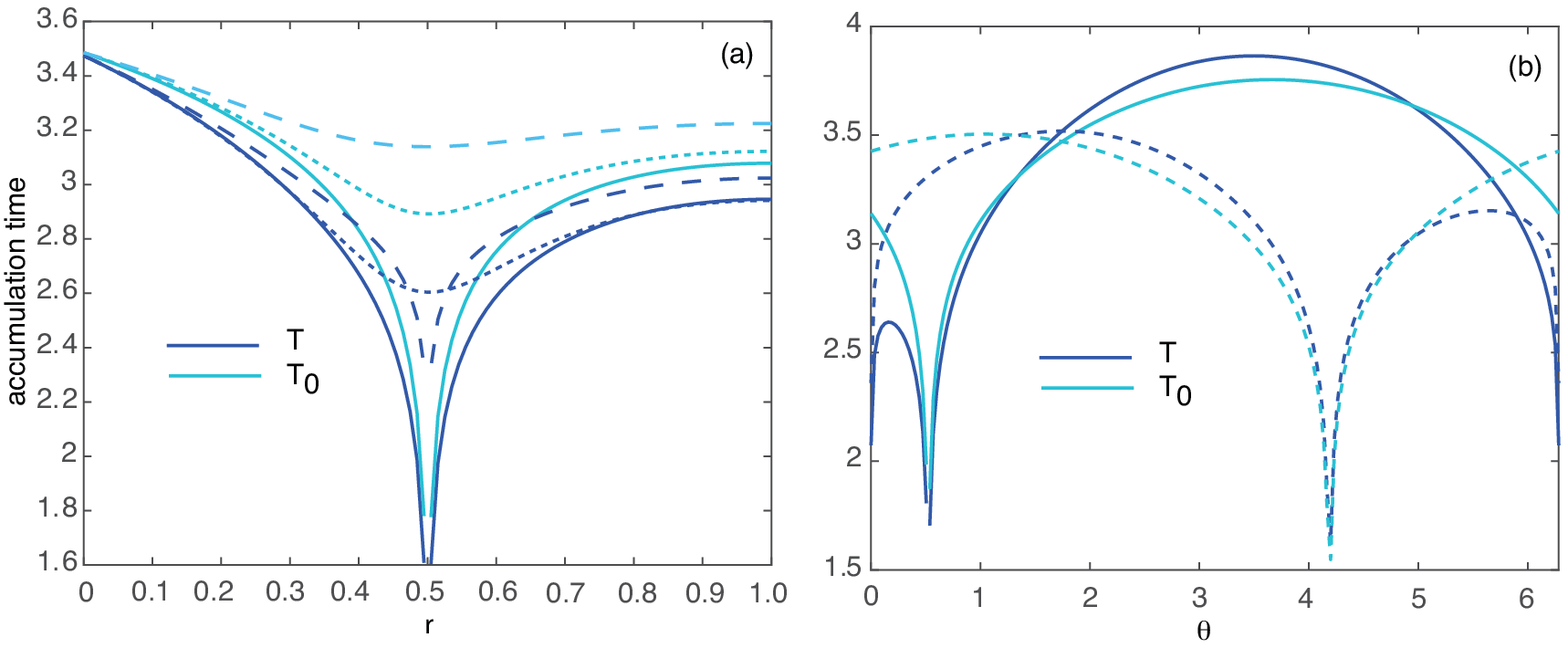} 
\caption{Accumulation time in the unit disc for the single-hole configuration shown in Fig. \ref{fig2}(a) with $\x_1=a(\cos \theta_1,\sin\theta_1)$ and $\x_0=(b,0)$. (a) Plots of $T(\x) $ and $T_0(\x)$ as a function of $r$ for $\x=r(\cos \theta,\sin \theta)$ and $\theta_1=\pi/6$ with $\theta=\pi/6$ (solid curves), $\theta=\pi/12$ (dotted curves) and $\theta=0$ (dashed curves). (b) Plots of $T(\x)$ and $T_0(\x)$ as a function of $\theta$ for $\x=a(\cos \theta,\sin \theta)$ with $\theta_1=\pi/6$ (solid curves) and $\theta_1=4\pi/3$ (dashed curves). Other parameter values are $\Gamma_0=1$, $\nu=0.1$, $a=b=0.5$ and $D=1$.}
\label{fig4}
\end{figure}

As our first example, consider the 2D configuration shown in Fig. \ref{fig2}(a). The domain $\Omega$ is taken to be the unit disc with a single hole placed at $\x_{1}=a(\cos \theta_1,\sin \theta_1)$. The boundary condition is
\begin{equation}
\nabla u(\x,t)\cdot {\bf n}=0,\ |\x|=1,\quad u(\x,t)=1,\ |\x-\x_{1}|=\epsilon.
\end{equation}
We have taken $\Phi_1=1$.
The initial concentration is localized at a point on the $x$-axis so that
\begin{equation}
u(\x,0)=\Gamma_0\delta(\x-\x_0),\quad \x_0=(b,0), \ 0 < b < 1, \Gamma_0 < \pi.
\end{equation}
In the case of the unit disc, the Neumann Green's function $G_0(\x,\xi)$ is known explicitly:
\begin{align}
G_0(\x,\hxi)&=\frac{1}{2\pi}\left [-\ln(|\x-\hxi|)-\ln\left (\left |\x|\hxi|-\frac{\hxi}{|\hxi |}\right |\right )  +\frac{1}{2}(|\x|^2+|\hxi|^2)-\frac{3}{4}\right ],
\label{Gdisc}
\end{align}
with the regular part obtained by dropping the first logarithmic term. It follows from equation (\ref{Tres1}) that the accumulation time is
 \begin{align}
T(\x)&= \frac{\pi-\Gamma_0}{2\pi \nu D}- \Gamma_0 \left [ G_0(\x,\x_0)-G_0(\x_1,\x_0)\right ]\\
&\quad -[\pi-\Gamma_0][G_0(\x,\x_1)-R_0(\x_1,\x_1)]+O(\nu).\nonumber 
\end{align}
In Fig. \ref{fig3}(a-c) we compare the full accumulation time $T(\x)$ with the truncated version $T_0(\x)$ for $a=0.5$, $b=0.5$ and $\theta_1=\pi/6$. It can be seen that the difference between the two is maximized in a neighborhood of the initial position. This is consistent with equation (\ref{Tres1}), which shows that the difference is proportional to $G_0(\x,\x_0)$ at $O(1)$. The full accumulation time in the case $a=0.5$, $b=0$ and $\theta_1=4\pi/3$ is shown in Fig. \ref{fig3}(d). In Fig. \ref{fig4}(a) we plot $T(\x)$ and $T_0(\x)$ as a function of $r$ for $\x=r(\cos\theta,\sin \theta)$ and $\theta$ fixed. Analogous plots of the accumulation times as a function of $\theta$ for fixed $r$ are shown in Fig. \ref{fig4}(b). These examples illustrate the singular nature of the accumulation time as $\x$ approaches the initial position $\x_0$ or the hole at $\x_1$.

\subsection{Pair of targets in the unit disc}

\begin{figure}[t!]
\centering
\includegraphics[width=14cm]{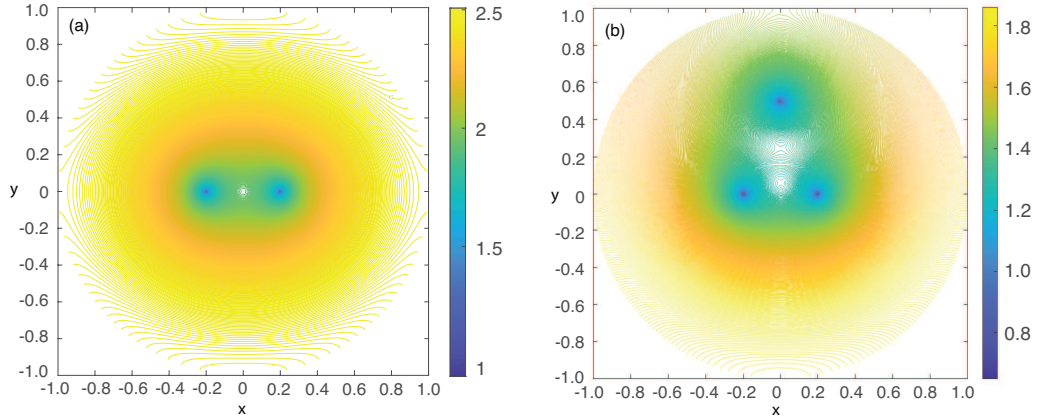} 
\caption{Contour plots of the accumulation time $T(\x)$ in the unit disc for the configuration shown in Fig. \ref{fig2}(b) with $\x_1=(a,0)$, $\x_2=(-a,0)$= and $\x_0=(0,b)$: (a) $\Gamma_0=0$ and (b) $\Gamma_0=1$. Parameter values are $\nu=0.1$, $D=1$, $a=0.2$ and $b=0.5$.}
\label{fig5}
\end{figure}

As our second example, consider the 2D configuration shown in Fig. \ref{fig2}(b). The domain $\Omega$ is again the unit disc, but now there is a pair of identical holes placed on the $x$-axis at $\x_{1}=(a,0)$ and $\x_2=(-a,0)$, $0<a<1$. The boundary conditions are
\begin{equation}
\nabla u(\x,t)\cdot {\bf n}=0,\ |\x|=1,\quad u(\x,t)=1,\ |\x-\x_{1,2}|=\epsilon.
\end{equation}
We have taken $\Phi_j=\overline{\Phi}=1$ for $j=1,2$.
The initial concentration is localized at a point on the $y$-axis so that
\begin{equation}
u(\x,0)=\Gamma_0\delta(\x-\x_0),\quad \x_0=(0,b), \ 0 < b < 1, \Gamma_0 < \pi.
\end{equation}
It follows from equation (\ref{one}) that the accumulation time for two identical targets in the unit disc with $\overline{\Phi}=1$ is given by
 \begin{align}
T(\x)&=\frac{1}{D}\frac{\pi-\Gamma_0}{4\pi \nu }- \Gamma_0 \left [ G_0(\x,\x_0)-\frac{1}{2}[G_0(\x_1,\x_0)+G_0(\x_2,\x_0)]\right ]\\
&\quad -\frac{\pi-\Gamma_0}{2}\bigg \{  G_0(\x,\x_1)+G_0(\x,\x_2)\nonumber \\
&\qquad \qquad -\frac{1}{2}\bigg [R_0(\x_1,\x_1)+R_0(\x_2,\x_2)+
G_0(\x_1,\x_2)+G_0(\x_2,\x_1)\bigg ]\bigg \}+O(\nu).\nonumber 
\end{align}
In Fig. \ref{fig5} we show contour plots of the $O(1)$ accumulation time $T(\x)$, $\x\in \Omega\backslash (\calU_1\cup \calU_2)$, for $a=0.2$ and $b=0.5$ in the two cases $\Gamma_0=0$ (zero initial concentration) and $\Gamma_0=1$. As expected, the plot is symmetric with respect to reflections about the $y$-axis. When $\Gamma_0=0$ there are local minima of $T(\x)$ in the vicinity of the holes, whereas there is an additional minimum at $\x_0$ when $\Gamma_0>0$.

\section{Discussion}
In this paper, we revisited the classical problem of analyzing the relaxation to steady state of diffusion in a 2D singularly perturbed domain. The usual approach is to calculate the principal non-zero eigenvalue of the negative Laplacian using matched asymptotics \cite{Ward93,Ward93a}. However, this does not allow for differences in the rate of relaxation at different spatial locations and loses all information about the initial position. Moreover, it relies on the existence of a sufficiently large spectral gap. Therefore, we developed an alternative method for characterizing the approach to steady-state based on the so-called accumulation time. Although one could consider an eigenfunction expansion of the accumulation time, which would yield a spatially varying relaxation rate, such an approximation still relies on a spectral gap. Here we showed how an asymptotic expansion of the full accumulation time can be obtained by solving the diffusion equation in Laplace space without any recourse to a spectral decomposition. The outer solution $\u(\x,s)$ was then used to calculate the accumulation time according to $T(\x)=-u^*(\x)^{-1}\lim_{s\rightarrow 0} d(s\u(\x,s))/ds$. We also highlighted several general features of $T(\x)$ . First, $T(\x)\rightarrow \infty$ as $\epsilon \rightarrow 0$. Second, for finite $\epsilon$, the accumulation time is more sensitive to the spatial locations of the holes than the steady-state concentration $u^*(\x)$. For example, terms depending on the spatial locations occur at $O(\nu)$ in the asymptotic expansion of $u^*(\x)$, whereas they occur at $O(1)$ in the case of $T(\x)$. Third, $T(\x)$ has a strong dependence on the initial condition that cannot be captured by focusing on the principal eigenvalue and eigenfunction. Although we illustrated the analysis using simple geometric configurations, one could consider more complicated geometries with multiple holes, provided that the associated Neumann Green's function and its regular part were known or could be determined numerically. 

There are a number of possible generalizations of the current work. First, within the context of 2D diffusion, one could consider more general exterior and interior boundary conditions, provided that there existed a unique steady-state solution. For example, modifying the exterior boundary condition would simply change the definition of the Green's function used in the outer solution. On the other hand, imposing constant flux conditions on the hole boundaries, say, would modify the inner solution and the corresponding singularity conditions for the outer solution. Second, one could consider non-spherical hole shapes, provided that the corresponding shape capacitances could be determined \cite{Ward93,Ward93a}. 

Finally, one could develop an analogous asymptotic analysis of the accumulation time for diffusion in 3D singularly perturbed domains. However, the details of the matched asymptotic analysis in 3D differs considerably from 2D, reflecting differences in the singular nature of the modified Helmholtz Green's function. That is, as $|\x-\x_0|\rightarrow 0$, 
\begin{equation}
G(\x,s|\x_0)\rightarrow -\frac{1}{2\pi D}\ln|\x-\x_0| \mbox{ for } d=2, \quad G(\x,s|\x_0)\rightarrow \frac{1}{4\pi D|\x-\x_0|} \mbox{ for } d=3.
\end{equation}
In 2D it was necessary to consider an asymptotic expansion in powers of $\nu=-1/\ln \epsilon$ at $O(1)$ in $\epsilon$, where $\epsilon$ specifies the relative size of the holes. However, taking the small-$s$ limit was relatively straightforward. On the other hand, performing an asymptotic expansion in powers of $\epsilon$ for 3D diffusion results in terms of order $O((\epsilon/s)^n)$, which are singular in the limit $s\rightarrow 0$ for fixed $\epsilon >0$. We have previously shown how to remove these singularities by considering a triple expansion of the solution in $\epsilon$, $s$ and $\Lambda\sim \epsilon /s$, and performing partial summations over geometric series in $\Lambda$ \cite{Bressloff21B}. Similar methods would be needed in order to determine the asymptotic expansion of the accumulation time in powers of $\epsilon$. The details will be presented elsewhere.

\end{document}